\documentstyle[twoside,11pt]{article}
  \newcommand{\bref}[5]{\noindent\parbox[t]{0.8cm}{#1}\parbox[t]{15.2cm}{{#2}{\it #3}{\bf #4}#5}\par}

\begin{document}
  \title{\bf {Optimal conditions for $(L_1;L_2)$ to be forcibly bigraphic}
  \thanks{This work is supported by the National Natural Science Foundation of
  China(NSFC11921001,11601108) and the National Key Research and Development Program of
  China(2018YFA0704701). The authors are grateful to professor C. Zong for his supervision and discussion. }}
\author{ $ {\mbox {Jiyun \,Guo}}$, $ {\mbox {Yuqin \,Zhang}}$ \\
{\small  College of Mathematics, Tianjin University,}
{\small Tianjin, 300072, P. R. China } \\
 {\small Email: 1020233004@tju.edu.cn.} }
\date{}
\maketitle
\begin{center}
\parbox{0.9\hsize}
{\small {\bf Abstract.}\ \ Let $L_1=([a_1,b_1],\ldots,[a_m,b_m])$
and $L_2=([c_1,d_1],\ldots,[c_n,d_n]$) be two sequences of intervals
consisting of nonnegative integers with $b_1\ge  \cdots\ge b_m$ and
$d_1\ge \cdots\ge d_n$. In this paper, we first give two optimal
conditions for the sequences of intervals $L_1$ and $L_2$ such that
each pair $(P;Q)$ with $P=(p_1,\ldots,p_m)$, $Q=(q_1,\ldots,q_n)$,
$a_i\le p_i\le b_i$ for $1\le i\le m$, $c_i\le q_i\le d_i$ for $1\le
i\le n$ and $\sum\limits_{i=1}^m p_i=\sum\limits_{i=1}^n q_i$ is
bigraphic. One of them is optimal sufficient condition and the other
one optimal necessary condition. We also present a characterization
of $(L_1;L_2)$ that is forcibly bigraphic on sequences of intervals.
This is an extension of the well-known theorem on bigraphic
sequences due to Gale and Ryser(see [4],[9]).


}
\end{center}

\section*{1. Introduction}
\hskip\parindent Throughout this paper, we consider only undirected
graphs without loops or parallel edges. Let $P_m=(p_1,\ldots,p_m)$
and $Q_n=(q_1,\ldots,q_n)$
 be two non-increasing sequences
of nonnegative integers. The pair $(P_m;Q_n)$ is said to be {\it
bigraphic} if there is a bipartite graph $G$ with partite sets
$X=\{x_1,\ldots,x_m\}$ and $Y=\{y_1,\ldots,y_n\}$ such that
$d_G(x_i)=a_i$ for $1\le i\le m$ and $d_G(y_i)=b_i$ for $1\le i\le
n$. In this case, $G$ is called a {\it realization} of $(P_m;Q_n)$.
Gale [4] and Ryser [9],independently, gave a  characterization for
$(P_m;Q_n)$ to be bigraphic.

{\bf Theorem 1.1.(Gale[4], Ryser[9])}. \ \  $(P_m;Q_n)$ is bigraphic
if and only if $\sum\limits_{i=1}^m p_i=\sum\limits_{i=1}^n q_i$ and
$$\sum_{i=1}^r q_i\le  \sum_{i=1}^m
min\{p_i,r\}\mbox{ for each $r$ with $1\le r\le n$.}\eqno{(1)}$$

 Analogous to the  case of bigraphic sequences, Tripathi et al. [3]
considered two sequences of intervals of nonnegative integers and
provided a corresponding characterization, which generalized Theorem
1.1.

Let $L_1=([a_1,b_1],\ldots,[a_m,b_m]$) and
$L_2=([c_1,d_1],\ldots,[c_n,d_n]$) be two sequences of intervals
consisting of nonnegative integers with $a_1\ge  \cdots\ge a_m$ and
$c_1\ge  \cdots\ge c_n$. We say that $(L_1;L_2)$ is bigraphic if
there exists a bipartite graph $G$ with partite sets
$X=\{x_1,\ldots,x_m\}$ and $Y=\{y_1,\ldots,y_n\}$ such that $a_i\le
d_G(x_i)\le b_i$ for $1\le i\le m$ and $c_i\le d_G(y_i)\le d_i$ for
$1\le i\le n$. In this case, $G$ is referred to as a realization of
$(L_1;L_2)$. Tripathi et al. [3] provided a characterization of
$(L_1;L_2)$ that is bigraphic.

{\bf Theorem 1.2.(Tripathi et al.[3])} . \ \ $(L_1;L_2)$ is
bigraphic
 if and only if

$$\sum_{i=1}^t c_i\le  \sum_{j=1}^m
min\{b_j,t\}\mbox{ for each $t$ with $1\le t\le n$}$$ and
$$\sum_{i=1}^s a_i\le  \sum_{j=1}^n min\{d_j,s\}\mbox{ for each $s$
with $1\le s\le m$}.$$

 The main object of the paper is to investigate
optimal necessary condition and optimal sufficient condition on two
sequences of intervals and give a characterization for $(L_1;L_2)$
to be forcibly bigraphic. The notion forcibly bigraphic pair will be
introduced later.

{\bf Theorem 1.3.} \ \ Let $L_1=([a_1,b_1],\ldots,[a_m,b_m]$) and
$L_2=([c_1,d_1],\ldots,[c_n,d_n]$) be two sequences of intervals
consisting of nonnegative integers with $b_1\ge  \cdots\ge b_m$ and
$d_1\ge  \cdots\ge d_n$. If the following two inequalities hold:
$$\sum_{i=1}^k b_i\le  \sum_{i=1}^n
min\{c_i,k\}\mbox{ for each $k$ with $1\le k\le m$,}\eqno{(2)}$$

$$\sum_{i=1}^l d_i\le  \sum_{i=1}^m
min\{a_i,l\}\mbox{ for each $l$ with $1\le l\le n$,}\eqno{(3)}$$
then each pair $(P;Q)$ with $P=(p_1,\ldots,p_m)$,
$Q=(q_1,\ldots,q_n)$, $a_i\le p_i\le b_i$ for $1\le i\le m$, $c_i\le
q_i\le d_i$ for $1\le i\le n$ and $\sum\limits_{i=1}^m
p_i=\sum\limits_{i=1}^n q_i$ is bigraphic.

It follows easily from Gale-Ryser Theorem that the sufficient
condition for each pair $(P;Q)$ being graphic is optimal.
Additionally, when it comes to the converse proposition of Theorem
1.3, we point out that the condition is sufficient but not
necessary, as can be seen by a counterexample.

{\bf Counterexample 1.} Take $L_1=([2,3],[1,2])$ and
$L_2=([1,2],[0,1])$, which satisfy $b_1\ge b_2$ and $d_1\ge d_2$. It
is easy to check that every$(P;Q)$ with $P=(p_1,p_2)$,
$Q=(q_1,q_2)$, $a_i\le p_i\le b_i$ for $1\le i\le 2$, $c_i\le q_i\le
d_i$ for $1\le i\le 2$ and $\sum\limits_{i=1}^m
p_i=\sum\limits_{i=1}^n q_i$ is bigraphic. However, (2) does not
hold for $k=1$ and $k=2$.

{\bf Theorem 1.4.} \ \ Let $L_1=([a_1,b_1],\ldots,[a_m,b_m]$) and
$L_2=([c_1,d_1],\ldots,[c_n,d_n]$) be two sequences of intervals
consisting of nonnegative integers with $b_1\ge  \cdots\ge b_m$ and
$d_1\ge \cdots\ge d_n$.
 If  each pair $(P;Q)$ with
$P=(p_1,\ldots,p_m)$, $Q=(q_1,\ldots,q_n)$, $a_i\le p_i\le b_i$ for
$1\le i\le m$, $c_i\le q_i\le d_i$ for $1\le i\le n$ and
$\sum\limits_{i=1}^m p_i=\sum\limits_{i=1}^n q_i$ is bigraphic, then
we have
$$\sum_{i=1}^k b_i\le  \sum_{i=1}^n
min\{c_i,k\}+|\sum\limits_{i=1}^m b_i-\sum\limits_{i=1}^n c_i|
\mbox{ for each $k$ with $1\le k\le m$},\eqno{(4)}$$
$$\sum_{i=1}^l d_i\le  \sum_{i=1}^m
min\{a_i,l\}+|\sum\limits_{i=1}^n d_i-\sum\limits_{i=1}^m a_i|
\mbox{ for each $l$ with $1\le l\le n$}.\eqno{(5)}$$

Again, it follows from Gale-Ryser Theorem that the necessary
condition for each pair $(P;Q)$ being graphic is optimal.
 Note that the condition is necessary but not
sufficient, as can be seen by  counterexample 2.

 {\bf Counterexample 2.} Take $L_1=([1,3],[2,3])$ and
$L_2=([1,2],[0,2])$. It is easy to check that (4) and (5) all hold.
However, $(P;Q)$ is not bigraphic with $P=(1,3)$ and $Q=(2,2)$.

Moreover, the above two theorems are generalizations of the results
in [5] due to Guo and Yin. In fact, there are a lot of articles on
the subject of lists of graphs (such as  Furuya and Yashima [1],
Goyal et al. [2], Lai and Hu[6],Tripathi et al. [10], Roberts [8]
and Rao [7] and so on).

Combining Theorem 1.3 with Theorem 1.4, we  arrive at the
characterization of $(L_1;L_2)$  that is forcibly bigraphic. Before
presenting it, we first give a definition.

{\bf Definition 1.5.} \ \ Let $L_1=([a_1,b_1],\ldots,[a_m,b_m]$) and
$L_2=([c_1,d_1],\ldots,[c_n,d_n]$) be two sequences of intervals
consisting of nonnegative integers with $b_1\ge  \cdots\ge b_m$ and
$d_1\ge \cdots\ge d_n$. $(L_1;L_2)$ is said to be forcibly bigraphic
if every pair $(P;Q)$ with $P=(p_1,\ldots,p_m)$,
$Q=(q_1,\ldots,q_n)$, $a_i\le p_i\le b_i$ for $1\le i\le m$, $c_i\le
q_i\le d_i$ for $1\le i\le n$ and $\sum\limits_{i=1}^m
p_i=\sum\limits_{i=1}^n q_i$ is bigraphic.

{\bf Theorem 1.6.} \ \ Let $L_1=([a_1,b_1],\ldots,[a_m,b_m]$) and
$L_2=([c_1,d_1],\ldots,[c_n,d_n]$) be two sequences of intervals
consisting of nonnegative integers with $b_1\ge  \cdots\ge b_m$,
$d_1\ge \cdots\ge d_n$, $\sum\limits_{i=1}^n d_i=
\sum\limits_{i=1}^m a_i$ and $\sum\limits_{i=1}^n c_i=
\sum\limits_{i=1}^m b_i$. Then $(L_1;L_2)$ is forcibly bigraphic if
and only if
$$\sum_{i=1}^k b_i\le  \sum_{i=1}^n
min\{c_i,k\}\mbox{ for each $k$ with $1\le k\le m$} $$ and
$$\sum_{i=1}^k d_i\le  \sum_{i=1}^m
min\{a_i,k\}\mbox{ for each $k$ with $1\le k\le n$.}$$

One can see that Theorem 1.6  generalizes  Gale-Ryser Theorem, which
corresponds to $a_i=b_i$ for $1\le i\le m$ and $c_j=d_j$ for $1\le
j\le n$.

 {\bf 2. Proof of Theorem 1.3} \ \

Choose an arbitrary pair $(P;Q)$ with  $P=(p_1,\ldots,p_m)$ and
$Q=(q_1,\ldots,q_n)$,
 $a_i\le p_i\le b_i$ for $1\le i\le m$ , $c_i\le q_i\le d_i$ for $1\le i\le n$ and
 $\sum\limits_{i=1}^m p_i=\sum\limits_{i=1}^n q_i$.
 Suppose (2) and (3) all hold, and now it suffices to show the pair $(P;Q)$ is bigraphic.

Rearrange the elements in $P=(p_1,\ldots,p_m)$ and
$Q=(q_1,\ldots,q_n)$ such that they are in nonincreasing order and
then denote $P'=(p'_1,\ldots,p'_m)$ and $Q'=(q'_1,\ldots,q'_n)$,
where $p'_1\ge \ldots \ge p'_m$ and $q'_1\ge \ldots \ge q'_n$. On
the other hand, the intervals in $L_1$ and $L_2$ are also rearranged
so that they are in coincidence with the above order, i.e.,
$L'_1=([a'_1,b'_1],\ldots,[a'_m,b'_m])$ and
$L'_2=([c'_1,d'_1],\ldots,[c'_m,d'_n])$, where $a'_i\le p'_i\le
b'_i$ for $1\le i\le m$ and $c'_i\le q'_i\le d'_i$ for $1\le i\le
n$, which together with (3) yield (6) below for each $r$ with $1\le
r\le n,$
$$
 \sum\limits_{i=1}^{r} q_i\le\sum\limits_{i=1}^{r} q'_i\le \sum\limits_{i=1}^{r} d'_i
 \le \sum\limits_{i=1}^{r} d_i
 \le \sum\limits_{i=1}^m min\{a_i,r\}
 = \sum\limits_{i=1}^m min\{a'_i,r\}
 \le \sum\limits_{i=1}^m min\{p'_i,r\}. \eqno{(6)}$$
Thus
 $\sum\limits_{i=1}^{r} q'_i
 \le \sum\limits_{i=1}^m min\{p'_i,r\} $
 for all $r$. From Theorem 1.1, we know  $(P';Q')$ is bigraphic, which is
equivalent to that $(P;Q)$ is bigraphic and so we are done. $\Box$ \\

{\bf 3. Proof of Theorem 1.4} \ \

For the proof, we first rearrange the  integers $a_1,\ldots,a_m$
such that they follow non-increasing order. Assume that $a'_1\ge
\ldots\ge a'_m$ and set $L'_1=([a'_1,b'_1],\ldots,[a'_m,b'_m])$,
then it is easy to see that $\sum\limits_{i=1}^m a'_i=
\sum\limits_{i=1}^m a_i $ and $ \sum\limits_{i=1}^m
min\{a'_i,r\}=\sum\limits_{i=1}^m min\{a_i,r\}$ for each $r$ with
$1\le r\le n$.

In the following, we just need to prove the inequality in (5), since
the idea used to verify the other one is analogous. Without loss of
generality, we may  assume that $n\ge m$, $b_i>a_i$ for all $i\le m,
d_i>c_i$ for all $i\le n$, $\sum\limits_{i=1}^n d_i\ge
\sum\limits_{i=1}^m a_i$ and set $\alpha= min\{a_1,\ldots,a_m\}.$

 If $\sum\limits_{i=1}^n d_i=
\sum\limits_{i=1}^m a_i,$ then let
$P=(p_1,\ldots,p_m)=(a'_1,\ldots,a'_m)$ and
$Q=(q_1,\ldots,q_n)=(d_1,\ldots,d_n)$. Clearly, $P$ and $Q$ are two
non-increasing sequences satisfying $a'_i\le p_i\le b'_i$ for $1\le
i\le m$, $c_i\le q_i\le d_i$ for $1\le i\le n$ and
$\sum\limits_{i=1}^m p_i=\sum\limits_{i=1}^n q_i$. By the assumption
of Theorem 1.4, $(P;Q)$ is bigraphic, and then it follows from (1)
that $\sum\limits_{i=1}^r d_i\le \sum\limits_{i=1}^m
min\{a'_i,r\}=\sum\limits_{i=1}^m min\{a_i,r\}$ for $1\le r\le n$.
Thus (5) holds.

 If
$\sum\limits_{i=1}^n d_i > \sum\limits_{i=1}^m a_i,$ set
$\sum\limits_{i=1}^n d_i - \sum\limits_{i=1}^m a_i =t,$ then it is
easy to see that $1\le t_2=t\le nd_1-m\alpha.$ Two cases now arise,
 depending on whether $1\le r<m$ or $m\le r\le n.$

{\bf Case 1.} Suppose $1\le r<m$. To prove (5), we have to divide
the situation into three subcases. For simplicity, we call $(P;Q)$ a
proper pair of $(L'_1;L_2)$ if $P$ and $Q$ satisfy the following
conditions: $a'_i\le p_i\le b'_i$ for $1\le i\le m$, $c_i\le q_i\le
d_i$ for $1\le i\le n$ and $\sum\limits_{i=1}^{m}
p_i=\sum\limits_{i=1}^{n} q_i$. Otherwise, we call $(P;Q)$ an
improper pair.

 {\bf Subcase 1.1.} If $1\le t\le m$, then let
$P=(p'_1,\ldots,p'_m)=(a'_1+1,\ldots,a'_t+1,a'_{t+1},\ldots,a'_m)$
and $Q=(q'_1,\ldots,q'_n)=(d_1,\ldots,d_n)$. Note that $p'_1\ge
\cdots\ge p'_m$, $q'_1\ge \cdots\ge q'_n$ and $(P;Q)$ is a proper
pair of $(L'_1;L_2)$. So by the assumption of Theorem 1.4, $(P;Q)$
is bigraphic and thus (1) follows, then we derive
$$\begin{array}{lll}
  \sum\limits_{i=1}^{r} d_i
  &=&\sum\limits_{i=1}^r q'_i\\
  &\le&\sum\limits_{i=1}^m min \{p'_i,r\}\\
  &=& \sum\limits_{i=1}^t min \{a'_i+1,r\}+\sum\limits_{i=t+1}^m min \{a'_i,r\}\\
  &\le& \sum\limits_{i=1}^t min \{a'_i,r\}+t+\sum\limits_{i=t+1}^m min \{a'_i,r\}\\
  &=& \sum\limits_{i=1}^m min \{a'_i,r\}+t\\
  &=& \sum\limits_{i=1}^m min \{a_i,r\}+t_2.
 \end{array}$$
Hence (5) holds for all $r<m$.

  {\bf  Subcase 1.2.}  If $m<t\le n,$ then $m\le m+n-t<n,$ which together
with $1\le r<m$ yield $r< m+n-t.$ Let
$P=(p'_1,\ldots,p'_m)=(a'_1+1,\ldots,a'_m+1)$ and
$Q=(q'_1,\ldots,q'_n)=(d_1,\ldots,d_r,\ldots,d_{n+m-t},d_{n+m-t+1}-1,\ldots,d_n-1)$.
It's easy to see that  $P$ and $Q$ are two non-increasing sequences
and $(P;Q)$ is a proper pair. Again, since $(P;Q)$ is bigraphic, (1)
holds and so for each $i\le n$,
$$\begin{array}{lll}
  \sum\limits_{i=1}^{r} d_i
  &=&\sum\limits_{i=1}^r q'_i\\
  &\le&\sum\limits_{i=1}^m min \{p'_i,r\}\\
  &=& \sum\limits_{i=1}^m min \{a'_i+1,r\}\\
  &\le& \sum\limits_{i=1}^m min \{a'_i,r\}+m\\
  &<& \sum\limits_{i=1}^m min \{a_i,r\}+t_2.
 \end{array}$$

  {\bf Subcase 1.3.} If $n<t\le nd_1-m\alpha$
  (obviously, if $nd_1-m\alpha\le n$, the proof is complete.),
let $P_1=(p'_1,\ldots,p'_m)=(a'_1,\ldots,a'_m)$ and
$Q_1=(q'_1,\ldots,q'_n)=(d_1-1,\ldots,d_n-1)$. If there exist  $l$,
$r$ and $w$ such that $d_i-1=0$ for all $i>l$, $d_i-1>0$
 for all $i\le l$ and $l-w+1=t-n$, where $1\le w<l $,
 then let  $Q_2=(q'_1,\ldots,q'_{w-1},q'_w-1,\ldots,q'_l-1,q'_{l+1},\ldots,q'_n)$
 $=(d_1-1,\ldots,d_{w-1}-1,d_w-2,\ldots,d_l-2,d_{l+1}-1,\ldots,d_n-1)=(d_1-1,\ldots,d_{w-1}-1,d_w-2,\ldots,d_l-2,0,\ldots,0).$
 Otherwise repeat the same steps until the deficiency $t$ is removed thoroughly. Then we may assume that
 $Q_j=(q^j_1,\ldots,q^j_n)=(d_1-(j-1),\ldots,d_{x-1}-(j-1),d_x-j,\ldots,d_y-j,0\ldots,0),$
 where $d_y-j$ is the last term of nonnegative integer in $Q_j$.
It can be seen that $P_1$ and $Q_j$ are non-increasing lists and $
\sum\limits_{i=1}^{m} p'_i= \sum\limits_{i=1}^{n} q^{j}_i$. If
$(P_1;Q_j)$ is improper, then the proof is complete.
Otherwise(namely, $c_i\le q^{j}_i \le d_i$ for $1\le i\le n$) we
proceed with our proof as follows. Since $(P_1;Q_j)$ is bigraphic,
(1) reduces to (7) for each $r$ with $1\le r\le n$,
$$ \sum\limits_{i=1}^{r} q^j_i \le \sum\limits_{i=1}^m min
\{p'_i,r\}.\eqno{(7)}$$ Now we need to divide the range of
$r(i.e.,1\le r\le n)$  into three parts.

{\bf  Subsubcase 1.3.1.} Suppose $1\le r<x (x\in (1,n])$, then by
(7), we have $\sum\limits_{i=1}^{r} [d_i-(j-1)]\le
\sum\limits_{i=1}^m min \{a'_i,r\}$, implying $\sum\limits_{i=1}^{r}
d_i\le \sum\limits_{i=1}^m min \{a_i,r\}+t_2$ for every $r$ with
$1\le r\le n$, where $r(j-1)\le t_2.$

{\bf  Subsubcase 1.3.2.} Suppose $x\le r\le y (y\in[1,n]),$ then it
follows from (7) that $[d_1-(j-1)]+\ldots +[d_{x-1}-(j-1)]+[d_x-j]+
\ldots +[d_r-j] \le \sum\limits_{i=1}^m min \{a'_i,r\},$ that is,
$\sum\limits_{i=1}^{r} d_i\le \sum\limits_{i=1}^m min \{a_i,r\}+t_2$
for each $r$ with $1\le r\le n$, where $rj-x+1\le t_2.$

{\bf  Subsubcase 1.3.3.} Suppose $y<r\le n$.  It is not difficult to
obtain that $ \sum\limits_{i=1}^{r} q^j_i =[d_1-(j-1)]+\ldots
+[d_{x-1}-(j-1)]+[d_x-j]+\ldots +[d_y-j]+0+\ldots
+0=\sum\limits_{i=1}^{y} d_i-yj+x-1 $. Denote $d_i=j_i$ for $y+1\le
i\le n$, then
$Q_j=(q^j_1,\ldots,q^j_n)=(d_1-(j-1),\ldots,d_{x-1}-(j-1),d_x-j,\ldots,d_y-j,0\ldots,0)=
(d_1-(j-1),\ldots,d_{x-1}-(j-1),d_x-j,\ldots,d_y-j,d_{y+1}-j_{y+1},\ldots,d_n-j_n).$
Note that $(j-1)(x-1)+j(y-x+1)+(j_{y+1}+\ldots
+j_n)=(yj+1-x)+(j_{y+1}+\ldots+j_n)=t_2$ and (7) gives
 $\sum\limits_{i=1}^{r} q^j_i=[d_1-(j-1)]+\ldots
+[d_{x-1}-(j-1)]+[d_x-j]+\ldots +[d_y-j]+[d_{y+1}-j_{y+1}]+\ldots
+[d_r-j_r] \le \sum\limits_{i=1}^m min \{a'_i,r\}.$  The latter is
equivalent to $\sum\limits_{i=1}^{r}
d_i-[(j-1)(x-1)+j(y-x+1)+(j_{y+1}+\ldots +j_r)]\le
\sum\limits_{i=1}^m min \{a'_i,r\},$  reducing to that
$\sum\limits_{i=1}^{r} d_i\le \sum\limits_{i=1}^m min
\{a_i,r\}+(yj-x+1)+(j_{y+1}+\ldots+j_r)\le \sum\limits_{i=1}^m min
\{a_i,r\}+t_2.$ Therefore, (5) holds for each $r$ with $1\le r<n$,
which means that (5) also holds for each $r$ with $1\le r<m$ since
$m\le n$. Now we like to point out that the constructive proof, in
the case, can be implemented as an algorithm to construct a suitable
sequence $Q_j.$

{\bf Case 2.} Suppose $m\le r\le n$. To prove (5), we have to
consider four subcases depending on the variation of $t$.

{\bf  Subcase 2.1.} If $1\le t<m,$ then it is easy to get (5) as the
case is similar to 1.1.

{\bf  Subcase 2.2.} If $m\le t\le m+n-r,$ then $r\le m+n-t.$  So we
may take $P=(a'_1+1,\ldots,a'_m+1)$ and
$Q=(d_1,\ldots,d_r,\ldots,d_{m+n-t},d_{m+n-t+1}-1,\ldots,d_n-1).$
Observe that $(P;Q)$ is proper and bigraphic. Then it is stemmed
from (1) that $
 \sum\limits_{i=1}^{r} d_i
 \le \sum\limits_{i=1}^m min\{a'_i+1,r\}
\le \sum\limits_{i=1}^m min \{a'_i,r\}+m
 \le\sum\limits_{i=1}^m
min \{a_i,r\}+t_2.
$

{\bf  Subcase 2.3.} If $m+n-r<t\le n,$ then $m+n-t<r\le n.$ Let
$P=(a'_1+1,\ldots,a'_m+1)$ and
$Q=(d_1,\ldots,d_{m+n-t},d_{m+n-t+1}-1,\ldots,d_r-1,\ldots,d_n-1).$
Evidently, $P$ and $Q$ are non-increasing and $(P;Q)$ is proper. By
the assumption of Theorem 1.4, $(P;Q)$ is bigraphic and then (1)
yields $\begin{array}{lll}
 \sum\limits_{i=1}^{m+n-t} d_i+ \sum\limits_{i=m+n-t+1}^{r} (d_i-1)
 &\le& \sum\limits_{i=1}^m min\{a'_i+1,r\},
\end{array}$
this is further transformed to $\begin{array}{lll}
 \sum\limits_{i=1}^{r} d_i-(r+t-m-n)
 &\le& \sum\limits_{i=1}^m min\{a'_i,r\}+m,
\end{array}$
so we derive that $$\begin{array}{lll}
 \sum\limits_{i=1}^{r} d_i
 &\le& \sum\limits_{i=1}^m min\{a'_i,r\}+r+t-n\le \sum\limits_{i=1}^{m}
 min\{a_i,r\}+t_2.
\end{array}$$
Hence (5) holds for any $r$ with $m\le r\le n$.

{\bf  Subcase 2.4.} If $n<t \le nd_1-m\alpha$( if $ nd_1-m\alpha\le
n$, the proof is complete.), then set
$P_1=(p'_1,\ldots,p'_m)=(a'_1,\ldots,a'_m)$ and
$Q_1=(q'_1,\ldots,q'_n)=(d_1-1,\ldots,d_n-1)$. Upon applying similar
arguments as 1.3, we can derive (5) and so the proof is complete.
 $\Box$


\vskip 0.1cm

\noindent{\bf References}

\bref{[1]}{M. Furuya and T. Yashima,  The existence of f-forests and
f-tree in graphs,} { Discrete Appl. Math.,}{ 254}{ (2019),
113--123.}

\vskip 0.1cm

 \bref{[2]}{D. Goyal, V. Jayapaul and V.Raman,
Elusiveness of finding degrees,} { Discrete Appl. Math.,}{ 286}{
(2020), 128--139.}

\vskip 0.1cm

\bref{[3]}{A. Garg, A. Goel and A. Tripathi,  Constructive
extensions of two results on graphic sequences,} { Discrete Appl.
Math.,}{ 159}{ (2011), 2170--2174.}

\vskip 0.1cm

\bref{[4]}{D. Gale, A theorem on flows in networks,} { Pac. J.
Math,}{ 7}{ (1957), 1073--1082.}

\vskip 0.1cm

\bref{[5]}{J. Y. Guo and J.H. Yin, A variant of Niessen's problem on
degree sequences of graphs,} { Discrete Math. Ther. Comp. Scie,}{
16}{ (2014), 287--292.}

\vskip 0.1cm

\bref{[6]}{C. H. Lai and L. L. Hu, Potentially $K_m-G$-graphical
sequences: a survey,} {  Czechoslovak Math. J,}{ 59}{ (2009),
1059--1075.}

\vskip 0.1cm

\bref{[7]}{S. B. Rao,  A survey of the theory of potentially
P-graphic and forcibly P-graphic degree sequences,} {Lecture Notes
in Math,No.} {855,Springer Verlag,}{ (1981), 417--440.}

\vskip 0.1cm

\bref{[8]}{A. Robert, Tree matching, } {J. Graph Theory
}{95}{(2020), 59--75.}

\vskip 0.1cm

\bref{[9]}{H. J. Ryser,  Combinatorial properties of matrices of
zeros and ones,} {Cana. J. Math,} { 9}{ (1957), 371--377.}

\vskip 0.1cm

\bref{[10]}{A. Tripathi and S. Vijay, A note on a theorem of
Erd$\ddot{o}$s and Gallai, } {Discrete. Math.,} { 265}{ (2003),
417--420.}

\end{document}